\documentclass[12pt]{article}
\usepackage{amsfonts,amsmath2000}
\def \beq {\begin{eqnarray}}
\def \eeq {\end{eqnarray}}
\def \beqn {\begin{eqnarray*}}
\def \eeqn {\end{eqnarray*}}

\newcommand{\halmos}{\rule{1ex}{1.4ex}}
\newcommand{\clos}[1]{{\rm clos} \,#1} 
\newcommand{\ital}[1]{{\it #1}}

\newcounter{for}[section]
\renewcommand{\theequation}{\thesection.\arabic{for}}

\newtheorem{itlemma}{Lemma}[section]
\newtheorem{itproposition}[itlemma]{Proposition}
\newtheorem{theorem}[itlemma]{Theorem}
\newtheorem{itcorollary}[itlemma]{Corollary}
\newtheorem{itremark}[itlemma]{Remark}
\newtheorem{itremarks}[itlemma]{Remarks}
\newtheorem{itdefinition}[itlemma]{Definition}
\newtheorem{itexample}[itlemma]{Example}

\newenvironment{fact}{\begin{itfact}\rm}{\end{itfact}}
\newenvironment{claim}{\begin{itclaim}\rm}{\end{itclaim}}
\newenvironment{lemma}{\begin{itlemma}}{\end{itlemma}}
\newenvironment{remark}{\begin{itremark}\rm}{\end{itremark}}
\newenvironment{remarks}{\begin{itremarks} \rm}{\end{itremarks}}
\newenvironment{corollary}{\begin{itcorollary}}{\end{itcorollary}}
\newenvironment{proposition}{\begin{itproposition}}{\end{itproposition}}
\newenvironment{definition}{\begin{itdefinition}\rm}{\end{itdefinition}}
\newenvironment{example}{\begin{itexample}\rm}{\end{itexample}}
\newenvironment{proof}{\noindent {\em Proof}.\ \
}{\hspace*{\fill}$\halmos$\medskip}
\newcommand{\be}[1]{\addtocounter{for}{1} \begin{equation}\label{#1}}
\newcommand{\ee}{\end{equation}}

\newcommand{\bl}[1]{\begin{lemma}\label{#1}}
\newcommand{\br}[1]{\begin{remark}\label{#1}}
\newcommand{\brs}[1]{\begin{remarks}\label{#1}}
\newcommand{\bt}[1]{\begin{theorem}\label{#1}}
\newcommand{\bd}[1]{\begin{definition}\label{#1}}
\newcommand{\bp}[1]{\begin{proposition}\label{#1}}
\newcommand{\bc}[1]{\begin{corollary}\label{#1}}
\newcommand{\bfact}[1]{\begin{fact}\label{#1}}
\newcommand{\bex}[1]{\begin{example}\label{#1}}

\newcommand{\ec}{\end{corollary}}
\newcommand{\efact}{\end{fact}}
\newcommand{\eex}{\end{example}}
\newcommand{\el}{\end{lemma}}
\newcommand{\er}{\end{remark}}
\newcommand{\ers}{\end{remarks}}
\newcommand{\et}{\end{theorem}}
\newcommand{\ed}{\end{definition}}
\newcommand{\ep}{\end{proposition}}
\newcommand{\epr}{\end{proof}}
\newcommand{\bpr}{\begin{proof}}
\newcommand{\bcl}[1]{\begin{claim}\label{#1}}
\newcommand{\ecl}{\end{claim}}

\newcommand{\ecs}{\end{corollary}}
\newcommand{\eers}{\end{exercise}}
\newcommand{\eexs}{\end{example}}
\newcommand{\eems}{\end{example}}
\newcommand{\els}{\end{lemma}}
\newcommand{\eles}{\end{lemmaex}}
\newcommand{\ets}{\end{theorem}}
\newcommand{\eds}{\end{definition}}
\newcommand{\eps}{\end{proposition}}

\newcommand{\bi}{\begin{itemize}}
\newcommand{\ei}{\end{itemize}}
\newcommand{\ben}{\begin{enumerate}}
\newcommand{\een}{\end{enumerate}}
\newcommand{\go}[1]
     {\raisebox{-1ex}{$\,\stackrel{\textstyle{\leadsto}}{\scriptstyle{#1}}\,$}}
\newcommand{\inter}[1]{\stackrel{\circ}{#1}}

\def\vbar{\mathchoice{\vrule height6.3ptdepth-.5ptwidth.8pt\kern-.8pt}
   {\vrule height6.3ptdepth-.5ptwidth.8pt\kern-.8pt}
   {\vrule height4.1ptdepth-.35ptwidth.6pt\kern-.6pt}
   {\vrule height3.1ptdepth-.25ptwidth.5pt\kern-.5pt}}
\def\fudge{\mathchoice{}{}{\mkern.5mu}{\mkern.8mu}}
\def\bbc#1#2{{\rm \mkern#2mu\vbar\mkern-#2mu#1}}
\def\bbb#1{{\rm I\mkern-3.5mu #1}}
\def\bba#1#2{{\rm #1\mkern-#2mu\fudge #1}}
\def\bb#1{{\count4=`#1 \advance\count4by-64 \ifcase\count4\or\bba A{11.5}\or
   \bbb B\or\bbc C{5}\or\bbb D\or\bbb E\or\bbb F \or\bbc G{5}\or\bbb H\or
   \bbb I\or\bbc J{3}\or\bbb K\or\bbb L \or\bbb M\or\bbb N\or\bbc O{5} \or
   \bbb P\or\bbc Q{5}\or\bbb R\or\bbc S{4.2}\or\bba T{10.5}\or\bbc U{5}\or
   \bba V{12}\or\bba W{16.5}\or\bba X{11}\or\bba Y{11.7}\or\bba Z{7.5}\fi}}

\def \qed {{\hspace*{\fill}$\halmos$\medskip}}
\def \Z {{\bb Z}}
\def \R {{\bb R}}
\def \C {{\cal{C}}}
\def \N {{\cal{N}}}
\def \u {{\bf u}}
\def \w {{\bf w}}
\def \ra {\rightarrow }
\def \under {\underline{a}}
\def \o {\omega}
\def \s {\sigma}
\def \T {{\cal{T}}}
\def \P{{\cal{P}}}
\def \LL {{\cal{L}}}
\def \Lp {{\cal{L}}^+}
\def \Lm {{\cal{L}}^-}
\def \Sp {S_t^+}

\def \M {{\cal{M}}}
\def \A {{\cal{A}}}
\def \H {{\cal{H}}}
\def \Te {{\cal{T}_{\epsilon}}}
\def \ss {S}
\def \Pm {P_{-}}
\def \Pmk {P_{-\kappa}}
\def \Pp {P_+}
\def \FF {{\cal{F}}}
\def \E {{\cal{E}}}
\def \TT {{\cal{T}}}
\def \D {D}
\def \dir {{\cal{D}}}
\def \dint {\int_0^T\!\!\int}
\def \dspace {\int\!\!\int\!\!}

\newcommand{\ba}[1]{\addtocounter{for}{1} \begin{eqnarray}\label{#1}}
\newcommand{\ea}{\end{eqnarray}}
\def\Box{{\hfill\hbox{\enspace${\square}$}} \smallskip}
\def\sqr#1#2{{\vcenter{\vbox{\hrule height .#2pt
                             \hbox{\vrule width .#2pt height#1pt \kern#1pt
                                   \vrule width .#2pt}
                             \hrule height .#2pt}}}}
\def\square{\mathchoice\sqr54\sqr54\sqr{4.1}3\sqr{3.5}3}
\def\pmb#1{\setbox0=\hbox{#1}%
   \kern-.025em\copy0\kern-\wd0
   \kern.05em\copy0\kern-\wd0
   \kern-.025em\raise.0433em\box0 }
\def\sqr#1#2{{\vcenter{\vbox{\hrule height.#2pt
     \hbox{\vrule width.#2pt height#1pt \kern#1pt
   \vrule width.#2pt}\hrule height.#2pt}}}}

\def\ZZ{{\mathbb Z}}   %
\def\tZZ{{\tilde\ZZ}}
\def\A{{\mathcal A}}
\def\B{{\mathcal B}}
\def\RR{{\mathbb R}}   %
\def\NN{{\mathbb N}}   %
\def\PP{{\mathbb P}}   %
\def\EE{{\mathbb E}}   %
\def\II{{\mathbb I}}   %

\def\e{\epsilon}                
\def\t{\tau_{\gamma}}
\def\e{\epsilon}
\def\d{\delta}
\def\l{\lambda}
\def\L{\Lambda}
\def\n{\nu_{\Lambda}}
\def\nur{\nu_{\rho}}
\def\mur{\mu_{\rho}}
\def\g{\gamma}
\def\G{\Gamma}
\def\a{\alpha}
\def\b{\beta}
\def\r{\rho}
\def\v{\varphi}
\def\p{\partial}
\def\m{\mu_N(\eta)}
\def\z{\zeta}
\def\bs{\backslash}
\def\K{\Sigma}

\begin{document}

\title{Regularity of quasi-stationary measures for simple exclusion
in dimension $d\ge 5$}
\author{Amine Asselah \\C.M.I., Universit\'e de Provence,\\
39 Rue Joliot-Curie, \\F-13453 Marseille cedex 13, France\\
asselah@gyptis.univ-mrs.fr
\\ and \\  Pablo A. Ferrari\\IME-USP,\\ P.B. 66281,
05315-970\\ S\~ao Paulo, SP, Brazil\\ pablo@ime.usp.br}
\date{}
\maketitle
\begin{abstract}
We consider the symmetric simple exclusion process on $\ZZ^d$,
for $d\ge 5$, and study the regularity of the quasi-stationary
measures of the dynamics conditioned on not occupying the origin.
For each $\rho\in ]0,1[$, we establish uniqueness of the density of 
quasi-stationary measures in $L^2(d\nur)$, 
where $\nur$ is the stationary measure of density $\rho$.
This, in turn, permits us to obtain sharp estimates
for $P_{\nur}(\tau>t)$, where $\tau$ is the first time the origin is occupied.
\end{abstract}

{\em Keywords and phrases}: quasi-stationary measures,
exchange processes, hitting time, Yaglom limit.

{\em AMS 2000 subject classification numbers}: 60K35, 82C22,
60J25.

{\em Running head}: Regularity of quasi-stationary measures

\section{Introduction}
Let $\{\eta_t:t\ge 0\}$ be the symmetric simple exclusion
process on $\ZZ^d$. In this process, there is at most one particle per site 
(i.e. the state space is $\Omega:=\{0,1\}^{\ZZ^d}$), 
and at rate one the contents of
neighboring sites are interchanged. The homogeneous Bernoulli product measures,
say $\nur$ with density $\rho\in[0,1]$, are invariant and reversible for this
process.  Let $\tau$ be the first time the origin of $\ZZ^d$ is occupied by a
particle.  We are interested in two issues: (i) to estimate the probability
that the origin remains empty for large time when the initial configurations
are drawn from $\nur$ for each $\rho\in ]0,1[$; (ii) to describe the law of
$\eta_t$, at large time $t$, conditioned on $\{\tau>t\}$, the event that the
origin is empty up to time $t$.

When the dimension of the lattice 
is larger than 4, we show that there exists a measure
$\mur$, such that for any continuous function $f$
\be{eq1}
\lim_{t\to\infty} E_{\nur}[f(\eta_t)|\tau >t]=\int f d\mur.
\ee
This establishes the so-called Yaglom limit \cite{yag}.
Such limiting measures can be intrinsically characterized as fixed points of
the semi-groups $\{T_t,t>0\}$ defined by
\be{eq2}
(T_t(\mu),f):=E_{\mu}[f(\eta_t)|\tau >t],\quad t>0.
\ee
Thus, fixed points of $\{T_t,\ t>0\}$ are dubbed quasi-stationary measures
\cite{SVJ}, \cite{fkm}.
Here, we study the regularity of $\mur$ and uniqueness when
the dimension $d>4$.
This, in turn, gives us sharp asymptotics for the probability of 
$P_{\nur}(\tau>t)$, namely
\be{eq3}
\lim_{t\to\infty} {\exp(-\lambda(\rho).t)\over P_{\nur}(\tau>t)}
=\int \left({d\mur\over d\nur}\right)^2 d\nur<\infty,
\ee
where $-\lambda(\rho)<0$ is the top of the spectrum in $L^2(\nur)$
of $\bar L$, the generator of the simple exclusion process
absorbed when hitting $\{\eta:\ \eta_0=1\}$.

We briefly summarize some relevant results of~\cite{ad}.
In dimensions 1 and 2, the
Yaglom limit is $\d_{\underline 0}$, the measure concentrated
on the configuration with no particle (and $\lambda(\rho)=0$).
In dimensions 3 and 4, $\lambda(\rho)>0$ for $\rho\in ]0,1[$,
and $\int_0^tT_s(\nur)ds/t$ converges to a quasi-stationary measure $\mur$. 
By analogy with the case of independent random walks~\cite{ad},
we conjecture that the Yaglom limit exists and that $\mur$ is singular
with respect to $\nur$. Thus, it is only for dimensions larger than 4
that we expect regularity of $\mur$ with respect to $\nur$.
\section{Notations and Results}
Henceforth, we consider dimensions larger or equal to 5, and $\rho\in
]0,1[$. The symmetric simple exclusion process (SSEP) on
the lattice $\ZZ^d$ can be graphically constructed ``\`a la Harris''~\cite{H} 
as follows. First, fix the initial configuration by assigning to each site of
$\ZZ^d$ a value in $\{0,1\}$ which indicates if the site is occupied or empty.
Then, to each bond --pairs of adjacent sites-- associate a Poisson processes
of intensity~1; Poisson processes of different bonds are independent and
independent of the initial configuration. At the times events (\emph{marks})
of each Poisson process, the values of the corresponding sites are interchanged.
In this way, each particle jumps when a mark is present; two particles may jump
at the same time in opposite directions. By labeling particles,
we can trace in time their trajectories: they evolve as the so-called
\emph{stirring particles}.  
This construction is described in Arratia \cite{ar1}. When
the labels of the stirring particles are disregarded one obtains only the
occupation numbers; in this case the resulting process, called $\eta_t$, has
infinitesimal generator
\[
Lf(\eta) = \sum_{i \in \ZZ^{d}} \sum_{j:j\sim i} [
f(\eta^{i,j}) - f(\eta)],\quad{\rm for}\quad \eta\in \{0,1\}^{\ZZ^{d}},
\]
where $\eta^{i,j}_k = \eta_k + (\d_{kj}- \d_{ki})(\eta_i-\eta_j)$
and $i\sim j$ means that $|i_1-j_1|+\dots+|i_d-j_d|=1$. It is
well known that the process is Feller, and the product measures of density
$\r$ in $[0,1]$, say $\nu_{\r}$, are reversible for $L$ (see Chapter VIII of
Liggett \cite{l1}). In other words, $L$ is an unbounded self-adjoint operator
in $L^2(d\nur)$, and local functions form a core for the domain, say $\D(L)$.
We denote by $P_{\nur}$ the law of the SSEP with initial measure $\nu_{\r}$.
Let $\A=\{\eta:\ \eta_0=1\}$ and denote by $\tau$ the time of first occurrence
of $\A$.  As we are interested in the Dirichlet problem on $\A^c$, we
introduce $\H_\A=\{ \v\in L^2(\nur):\ \v(\eta)=0$ for $ \eta\in \A\}$.  Let
$\bar L$ be the operator defined by
\[
\bar Lf=1_{\A^c}Lf,\quad{\rm for}\quad f\in \D(L)\cap \H_\A.
\]
This corresponds to the simple exclusion dynamics absorbed when hitting the
event $A$. 
$\bar L$ is self-adjoint on $\H_\A$ with respect to $\nur$. 
We call $\{\bar S_t, t>0\}$ the
corresponding sub-Markovian semi-group of bounded 
operators on $L^2(\H_\A,\nur)$. In other words,
\[
\forall t>0,\quad \bar S_tf(\eta)=E^{\eta}[f(\eta_t)1_{\{\tau>t\}}].
\]
We denote by $T_t(\nur)$ the probability measure defined by duality
on $\v\in \H_\A$, 
\[
(T_t(\nur),\v)={\int \bar S_t\v\, 1_{\A^c}\,d\nur\over P_{\nur}(\tau>t)}=
\int \v {\bar S_t1_{\A^c}\over P_{\nur}(\tau>t)}d\nur.
\]
where we have used reversibility to obtain the third term.
Thus, if $f_t$ is the density of $T_t(\nur)$ with respect to $\nur$
\be{eq0.1}
\forall t>0, \quad f_t(\eta)=
{\bar S_t1_{\A^c}(\eta)\over P_{\nur}(\tau>t)}=
{P^{\eta}(\tau>t)\over \int\! P^{\zeta}(\tau>t)d\nu_{\rho}(\zeta)}.
\ee
It was established in \cite{ad} that a non-trivial quasi-stationary measure,
say $\mur$, could be obtained as limit along a Ces\`aro subsequence of
$T_t(\nur)$. Our main result is the following theorem.  
\bt{the1} If the
dimension is larger or equal to 5, then $\mur$ is absolutely continuous with
respect to $\nur$.  Moreover, for any integer $k\ge 1$, $d\mur/d\nur\in
L^k(\nur)$ and
\be{main}
\lim_{t\to\infty}\int\left({dT_t(\nur)\over d\nur}-
{d\mur\over d\nur}\right)^2d\nur=0.
\ee
\et
\br{rem1} This is stronger than establishing the Yaglom limit,
i.e. $\lim T_t(\nur)=\mur$.
As a consequence, $f:=d\mur/d\nur$ belongs to $\D(\bar L)$ and
satisfies (in the $L^2(\nur)$-sense)
\be{eq0.4}
\bar Lf+\lambda(\rho) f=0,\quad {\rm and}\quad \bar S_t f=e^{-\lambda(\rho) t} f,
\ee
with (see Theorem 2 of \cite{ad})
\be{bis-e}
\lambda(\rho)= \inf\Bigl\{
{(f, -\bar L f)_{\nur}\over (f,f)_{\nur}}:\ f\in\D(L)\cap \H_\A\Bigr\}.
\ee
\er
Theorem~\ref{the1} is based on the apriori bounds
through the
following lemma, in which we reformulate a general result
essentially contained in~\cite{ad}.
\bl{lem.gene} Let $\{\bar S_t\}$ be the semi-group
of a process absorbed when hitting a set $\A$.
Assume that $\{\bar S_t\}$ is reversible with respect
to $\nu$, and let $\lambda<\infty$ be given by \eqref{bis-e}.
The following two conditions are equivalent
\be{eq.gene}
(i)\ \sup_{t>0} \frac{e^{-\lambda t}}{P_{\nu}(\tau>t)}<\infty,
\quad{\rm and}\quad
(ii)\ \sup_{t>0} \int f_t^2d\nu<\infty.
\ee
Moreover, if either $(i)$ or $(ii)$ holds, then the Yaglom limit
$\mu$ exists and \eqref{main} holds.
\el
Now, the apriori bounds is a corollary of the following proposition, interesting
on its own.
\bp{prop1}Let the dimension $d\ge 3$.
Let $i\in \ZZ^d\bs\{0\}$ and $\eta\in \Omega$ with
$\eta_i=0$. We denote by $\eta^i$ the configuration identical
to $\eta$ except on $i$, where its value is 1. 
There is a constant $C_d$, independent of $i$ and $\eta$ such that
for any $t>0$,
\be{keypoint}
0\le P^{\eta}(\tau>t)-P^{\eta^i}(\tau>t)\le
C_dP^{\eta^i}(\tau>t)\PP(H_i<\infty),
\ee
where $H_i$ denotes the first time a symmetric random walk
starting at $i$ hits the origin.
\ep
The relation (\ref{keypoint}) would be obvious if the particles
were independent. Though it is rather intuitive for the symmetric
exclusion, our proof is rather long. A sketch of it is as follows.
We first write $P^{\eta}(\tau>t)$ in terms of a dual process, 
say $\{ X(\emptyset,t)\}$, which
corresponds to a stirring process on $\ZZ^d\bs\{0\}$ with birth at 
the nearest neighbors of the origin and with
initial condition an empty configuration. Then,
$P^{\eta}(\tau>t)-P^{\eta^i}(\tau>t)$ corresponds to the weight of
all paths whose end-points $X(\emptyset,t)=U\cup \{i\}$ with $\eta_j=0$ for
all $j\in U\subset \ZZ^d\bs\{0,i\}$. The problem is then to uncouple
$U$ from $\{i\}$. We then re-express $P(X(\emptyset,t)=U\cup \{i\})$ in terms of
a dual {\bf with finitely many particles}, say $\{\L_t\}$. Note that
$\{\L_t\}$ is not the `natural dual' of $\{ X(\emptyset,t)\}$ and
the correspondence is obtained through a
Feynman-Kac formula. Then we show a correlation inequality for the
expression in terms of
$\{\L_t\}$ by generalizing Andjel's inequality~\cite{an}.

The results about apriori bounds is the following.
\bc{lem1} Let the dimension $d\ge 5$. (i) There is a product measure
$\nu_{\a(.)}$ of density $\a(i)$ for $i\in \ZZ^d$ such that for any $t>0$
\[
\nu_{\a(.)}\prec T_t(\nur) \prec \nur,\quad{\rm and}\quad
\sum_{i\in \ZZ^d}(1-\frac{\a_i}{\rho})^2<\infty, 
\]
where $\prec$ denotes stochastic domination.

(ii) For any integer $k\ge 1$, 
there is a positive constant $C$, such that
\be{eq0.3}
\sup_{t>0} \int f_t^k(\eta) d\nu_{\rho}(\eta)\le C.
\ee
\ec
A consequence of Lemma~\ref{lem.gene} 
is a sharp asymptotic estimate for the tail of $\tau$ 
(compare with \cite{ad} Lemma 1).  
\bc{cor1} If the dimension $d\ge 5$, then
\be{eq0.5} \lim_{t\to\infty} {e^{-\lambda(\rho) t}\over
P_{\nur}(\tau>t)}=\int f^2 d\nur.  
\ee 
\ec 
Finally, we have a uniqueness result and some properties of $\mur$.
\bt{the2}(i) 
The set $\{\mu\ll \nur\,:\, \mu$ is quasi stationary and
$d\mu/d\nur\in L^2(\nur\}$ contains only $\mur$. 
(ii) For $i\in \ZZ^d$, define $\theta_i:\Omega\to\Omega$ with
$\theta_i\eta_k=\eta_{k+i}$, then
for any $\v$ local (i.e. depending on finitely many sites)
\[
\lim_{||i||\to\infty} \int \v(\theta_i\eta) d\mur(\eta)=\int \v d\nur.
\]
(iii) If $\nu$ is a probability with a continuous density with
respect to $\nur$, then
\[
\lim_{t\to\infty}{1\over t}\int_0^t T_s(\nu)ds=\mur.
\]
The convergence holds in weak-$L^2(\nur)$.
\et
Proposition~\ref{prop1} is proven in Section~\ref{regular}. 
Section~\ref{apriori} contains the proofs of 
Corollary~\ref{lem1}, and of Theorem~\ref{the1}.
In section~\ref{unique}, we establish the uniqueness 
part of Theorem~\ref{the2}.
In section~\ref{density}, we show that in $\mur$ the density at infinity
is $\rho$, and we conclude with the result about the basin of attraction
of $\mur$. 
\section{Proof of Proposition~\ref{prop1}.}
\label{regular}
\subsection{Duality and Feynman-Kac.}
We first express $P^{\eta}(\tau>t)$ using the \emph{dual} process (\cite{l1},
\cite{ar1}) based on the fact that the Poisson clocks associated to bonds are
invariant by time reflections. The dual process tracing back-in-time the
positions of the stirring particles touching the origin can be described using
the graphical construction at the beginning of Section 2. Again at each bond,
there is an independent mark process corresponding to the realization of a
Poisson process of intensity~1. At each mark between 0 and one of its nearest
neighbor, say $i$, a particle is born at $i$ unless $i$ is already occupied
(in which case nothing happens); the particles born in this way evolve
afterwards as stirring particles on $\ZZ^d\bs\{0\}$.  The only difference with
the previous construction is that now it is imposed to the origin to be
occupied at all times ---so that when it becomes empty, it is immediately
occupied with a newly created particle. Assume that at time 0, the lattice is
empty and let $X(\emptyset,t)$ be the set of sites 
occupied by the stirring particles at
time t; all these particles have been created at the origin. Let $\PP$ denote
averages over the Poisson realizations.  If $\P^*$ is the collection of finite
subsets of $\ZZ^d\bs \{0\}$, then the duality formula reads for any $t>0$
\be{eq1.1}
P^{\eta}(\tau>t)\;=\;(1-\eta_0)\sum_{\L\in \P^*}
\PP\left(X(\emptyset,t)=\Lambda\right)\prod_{j\in \L} (1-\eta_j).  
\ee 
Thus, if $\eta$ is such that $\eta_i=0$
\[
P^{\eta}(\tau>t)-P^{\eta^i}(\tau>t)\;=\;(1-\eta_0)\sum_{\L\in \P^*,\ \L\ni i}
\PP\left(X(\emptyset,t)=\Lambda\right)\prod_{j\in \L} (1-\eta_j).
\]
Assume, for a moment, that for $U\in\P^*$ and $i\not\in U$, we have
a constant $C_d$ independent of $i$ and $t$ such that
\be{new.ineq}
\PP\left(X(\emptyset,t)=U\cup\{i\}\right)\le C_d
\PP\left(X(\emptyset,t)=U\right)
\PP(H_i<\infty),
\ee
where we denote by $H_i$ the first time a symmetric random walk 
starting at $i$ hits the origin. Then, for $\eta$ such that $\eta_i=0$
\[
P^{\eta}(\tau>t)-P^{\eta^i}(\tau>t)\;\le\;C_d\PP(H_i<\infty)
(1-\eta_0)\sum_{U\in \P^*, i\not\in U}
\PP\left(X(\emptyset,t)=U\right)\prod_{j\in U\cup\{i\}} (1-\eta_j) 
\]
\[
\quad\quad\le C_d(1-\eta_i)\PP(H_i<\infty)P^{\eta^i}(\tau>t).
\]
Thus, it remains to prove (\ref{new.ineq}).

Let $\Lp$ be the generator of $\{X(\emptyset,t),\ t\ge 0\}$, and
let $\Sp$ be the associated semi-group. We first express the
dual of $\{X(\emptyset,t),\ t\ge 0\}$ in terms of a process
with finitely many particles. Actually, we are only interested
in $\Sp(1_{\L})(\emptyset):=\PP(X(\emptyset,t)=\L)$ for $\L\in\P^*$.
Let $\L$ and $A$ be in $\P^*$. We have, using $\Delta$ for the symmetric
difference,
\[
\Lp(1_{\L})(A)=\sum_{
\substack{\scriptstyle{x\sim y;\ x,y\not=\{0\}}\\
\scriptstyle{|A\Delta\{x,y\}|=|A|}}}
[1_{\L}(A\Delta\{x,y\})-1_{\L}(A)]+\sum_{
\substack{\scriptstyle{y\sim 0}\\ \scriptstyle{y\not\in A}}}
[1_{\L}(A\cup\{y\})-1_{\L}(A)].
\]
The first sum corresponds to the stirring process over $\ZZ^d\bs\{0\}$,
while the second sum corresponds to birth at the origin.
We reexpress now the last sum. For simplicity, we omit to write $y\sim 0$.
Thus,
\[
\sum_{y\not\in A}[1_{\L}(A\cup\{y\})-1_{\L}(A)]
=\sum_{\substack{\scriptstyle{y\not\in A}\\
\scriptstyle{y\in \L}}}[1_{\L}(A\cup\{y\})-1_{\L}(A)]-
\sum_{\substack{\scriptstyle{y\not\in A}\\
\scriptstyle{y\not\in \L}}}1_A(\L).
\]
We claim that this expression is equal to
\[
\C:=\sum_{y\in \L}[1_{A}(\L\bs\{y\})-1_{A}(\L)]-
\sum_{y\not\in \L}1_{A}(\L)+\sum_{y\in \L}1_{A}(\L).
\]
Indeed, we expand $\C$
\beqn
\C&=&
\sum_{\substack{\scriptstyle{y\in \L}\\
\scriptstyle{y\not\in A}}}[1_{A}(\L\bs\{y\})-1_{A}(\L)]+
\sum_{\substack{\scriptstyle{y\in \L}\\
\scriptstyle{y\in A}}}[1_{A}(\L\bs\{y\})-1_{A}(\L)]
-\sum_{y\not\in \L}1_{A}(\L)+\sum_{y\in \L}1_{A}(\L)\cr
&=& \sum_{\substack{\scriptstyle{y\not\in A}\\
\scriptstyle{y\in \L}}}[1_{A}(\L\bs\{y\})-1_{A}(\L)]
-\sum_{y\not\in \L}1_{A}(\L).
\eeqn
Thus, calling $\N_0:=\{y\in\ZZ^d: y\sim 0\}$ and using the self-duality
of the stirring part of $\Lp$, we obtain
\be{def-dual}
\Lp(1_{\L})(A)=\sum_{
\substack{\scriptstyle{x\sim y;\ x,y\not=0}\\
\scriptstyle{|\L\Delta\{x,y\}|=|\L|}}}
[1_{A}(\L\Delta\{x,y\})-1_{A}(\L)]
+ \sum_{\substack{\scriptstyle{y\sim 0}\\
\scriptstyle{y\in \L}}}[1_{A}(\L\bs\{y\})-1_{A}(\L)]+V(\L)1_A(\L).
\ee
where we set $V(\L):=2|\N_0\cap \Lambda|-|\N_0|$. Now, let
$\Lm$ denote the generator of the stirring process on $\ZZ^d\bs\{0\}$
with death when particles jump on the origin. Then, (\ref{def-dual})
can be written like
\[
\Lp(1_{\L})(A)=\Lm(1_A)(\L)+V(\L)1_A(\L).
\]
Thus, if $u(\L,t):=\Sp(1_{\L})(\emptyset)$, we have
\[
{du(\L,t)\over dt}=\Sp(\Lp 1_{\L})(\emptyset)=\Lm u(\L,t)+V(\L)u(\L,t).
\]
Now, we call $\{\L(t), t\ge 0\}$ the process generated by $\Lm$ and
we use Feynman-Kac to obtain
\[
\PP(X(\emptyset,t)=\L)=
\EE^{\L}[e^{\int_0^t V(\L(s))ds} 1_{\emptyset}(\L(t))]=
e^{-|\N_0|t}\EE^{\L}[\exp(2\int_0^t |\N_0\cap \L(s)|ds)
1_{\emptyset}(\L(t))].
\]
Let now $U\in\P^*$ and $i\not\in U$. We show in section~\ref{correlation}
that if $\L\in \P^*$ and
\be{enrique-type}
g(\L,t):=\EE^{\L}[\exp(2\int_0^t |\N_0\cap \L(s)|ds)
1_{\emptyset}(\L(t))],\quad{\rm then}\quad
g(U\cup\{i\},t)\le g(U,t)g(\{i\},t).
\ee
Then, in section~\ref{upperbound} we prove that $g(\{i\},t)
\le C_d P(H_i<\infty)$ for $d\ge 3$. Inequality (\ref{new.ineq}) 
follows then readily.
\subsection{A generalized correlation inequality.}
\label{correlation}
To make the notations closer to those of Andjel~\cite{an},
we set $p(x,y)=1$ when $x\sim y$, and $p(x,y)=0$ otherwise. Also,
we realize our stirring process as an exclusion process: the
particles attempt to jump to one of their nearest neighboring sites
at the time marks of independent Poisson processes of intensity $2d$; if the
site chosen (each neighboring site is chosen with the same probability)
for the attempt is occupied, the particle stays still.
As we are blind to the labeling of particles, the trajectories
are, in law, indistinguishable from our initial stirring process.

We proceed by induction on $n$ to prove
that for any sets $A,B\in\P^*$, with 
$A\cap B=\emptyset$, for any $t>0$, $\a\in \RR$, and any $n$-tuples
$0\le s^1<s^2<\dots<s^n\le t$, the following inequality holds
\be{ineq.land}
\EE^{A\cup B}[e^{\a\sum_{k=1}^n |\N_0\cap \L(s^k)|}
1_{\emptyset}(\L(t))]\le
\EE^{A}[e^{\a\sum_{k=1}^n |\N_0\cap \L(s^k)|}
1_{\emptyset}(\L(t))]\EE^{B}[e^{\a\sum_{k=1}^n |\N_0\cap \L(s^k)|}
1_{\emptyset}(\L(t))].
\ee
This will be our induction hypothesis at order $n$.
Once (\ref{ineq.land}) is proved, inequality (\ref{enrique-type}) follows
easily as in Proposition 4.1 of ~\cite{landim}.
 
\noindent{\bf Step n=0.} 
We need to prove that for $A,B\in \P^*$ with $A\cap B=\emptyset$
\be{eq1.7}
\PP^{A\cup B}(\L(t)=\emptyset)\le \PP^A(\L(t)=\emptyset)
\PP^B(\L(t)=\emptyset).
\ee
Following an idea of Arratia~\cite{ar}, we represent the process $\L(t)$ as
limit of a stirring process with no absorption on an enlarged lattice: we link
the origin $0$ with $\tilde 0$, the origin of a three dimensional lattice
$\tZZ^3$ isomorphic to $\ZZ^3$ (here we fix $\ZZ^3$ for concreteness; any
graph supporting the stirring construction, for which the corresponding random
walk is transient would fit). On each bond of $\tZZ^3$ and on the bond
$(0,\tilde 0)$, the rates of stirring are set equal to $\kappa$ large.  On the
enlarged lattice $\ZZ^d\cup \tZZ^3$, the particles perform a conservative
stirring process, though with different rates whether they jump across the
bonds of $\ZZ^d$ or across the bonds of $\tZZ^3$ and $(0,\tilde0)$.  When a
particle hits the origin 0, it has a probability going to 1 as
$\kappa\to\infty$ to wander in $\tZZ^3$ up to time $t$ without using bonds of
$\ZZ^d$. We call $U(t)$ the stirring process on $\ZZ^d\cup\tZZ^3$, and
$\PP_{\kappa}$ the law of the Poisson marks on the enlarged lattice.
It is not difficult to show that for any $\Lambda\in \P^*$, 
\be{eq1.4}
\PP^{\L}\left(\L(t)=\emptyset\right) =\lim_{\kappa\to\infty} 
\PP_{\kappa}^{\L}(U(t)\subset \tZZ^3).
\ee
Now, for the stirring process on the enlarged lattice $\ZZ^d\cup\tZZ^3$,
we use a correlation inequality due to Andjel~\cite{an}: 
\be{eq1.5}
\PP_{\kappa}^{A\cup B}(U(t)\subset \tZZ^3)\;\le\; 
\PP_{\kappa}^A(U(t)\subset \tZZ^3)\,\PP_{\kappa}^B(U(t)\subset \tZZ^3).
\ee
Thus, after taking the limit $\kappa\to\infty$, we obtain (\ref{eq1.7}).

\noindent{\bf Step n}. Our proof follows essentially
Andjel's proof. Our induction hypothesis is that
(\ref{ineq.land}) is valid for $n-1$ instants of time.
Let $0\le s^1<s^2<\dots<s^n\le t$ be $n$ instants of time,
and for each $\L\in \P^*$ let
\[
g_n(\L,t;s^1,\dots,s^n)=\EE^{\L}[
\exp(\a \sum_{i=1}^n |\N_0\cap \L(s^i)|) 1_{\emptyset}
(\L(t))].
\]
We set $\l=2d(|A|+|B|)$ and we let $\tau_1$ be the first time
a particle of $A\cup B$ attempts a jump (i.e. $\tau_1$ is
an exponential time of parameter $\lambda$). Note that by the
Markov property
\begin{align}
\EE^{A\cup B}[&\exp(\a \sum_{i=1}^n |\N_0\cap \L(s^i)|) 1_{\emptyset}
(\L(t))1_{\{\tau_1>s^1\}}]\notag\\
&= P(\tau_1>s^1)e^{\a(|A\cap \N_0|+|B\cap \N_0|)}\EE^{A\cup B}[
\exp(\a \sum_{i=2}^n |\N_0\cap \L(s^i-s^1)|) 1_{\emptyset}
(\L(t-s^1))]\notag\\
&=P(\tau_1>s^1)e^{\a(|A\cap \N_0|+|B\cap \N_0|)}g_{n-1}(A\cup B,t-s^1;
s^2-s^1,\dots,s^n-s^1).\notag
\end{align}
Following~\cite{an}, using the shorthand notation
${\bf s}$ for $s^1,\dots,s^n$
(and its abuse ${\bf s-u}=(s^1-u,\dots,s^n-u)$), 
and writing $A\tilde \Delta\{x,y\}$ to mean $A\Delta\{x,y\}\bs\{0\}$,
for we have to account for deaths of particules when they jump on 0,
we have 
\begin{align}
g_n&(A\cup B,t;{\bf s})=P(\tau_1>s^1)
e^{\a|A\cap \N_0|}e^{\a|B\cap \N_0|}g_{n-1}(A\cup B,t-s^1;
s^2-s^1,\dots,s^n-s^1)\notag\\
&+\!\!\int_0^{s^1}du \l e^{-\l u} \frac{1}{\l}\Big\{
[\sum_{x,y\in A} p(x,y)+\sum_{x,y\in B} p(x,y)]
g_{n}(A\cup B,t-u;{\bf s-u})\notag\\
&+[\sum_{\substack{x\in A\\y\not\in A\cup B}}p(x,y)
g_{n}(A\tilde \Delta\{x,y\}\cup B,t-u;{\bf s-u})]\notag\\
&+[\sum_{\substack{x\in B\\y\not\in A\cup B}}p(x,y)
g_{n}(A\cup B\tilde \Delta\{x,y\},t-u;{\bf s-u})]\notag\\
&+[\sum_{\substack{x\in A\\y\in B}}p(x,y)+
\sum_{\substack{x\in B\\y\in A}}p(x,y)]
g_{n}(A\cup B,t-u;{\bf s-u})\Big\}\notag
\end{align}
Reasoning as if the particles in $A$ 
were independent from the particles in $B$, we obtain 
\begin{align}
g_n&(A,t;{\bf s}) g_n(B,t;{\bf s})=\notag\\
&P(\tau_1>s_1)
e^{\a|A\cap \N_0|}e^{\a|B\cap \N_0|}g_{n-1}(A,t-s^1;
s^2-s^1,\dots)g_{n-1}(B,t-s^1;s^2-s^1,\dots)\notag\\
&+\int_0^{s^1}du \l e^{-\l u} \frac{1}{\l}\Big\{
[\sum_{x,y\in A} p(x,y)+\sum_{x,y\in B} p(x,y)]
g_{n}(A,t-u;{\bf s-u})g_{n}(B,t-u;{\bf s-u})\notag\\
&+[\sum_{\substack{x\in A\\y\not\in A\cup B}}p(x,y)
g_{n}(A\tilde \Delta\{x,y\},t-u;{\bf s-u})g_{n}(B,t-u;{\bf s-u})]\notag\\
&+[\sum_{\substack{x\in B\\y\not\in A\cup B}}p(x,y)
g_{n}(A,t-u;{\bf s-u})g_{n}(B\tilde \Delta\{x,y\},t-u;{\bf s-u})]\notag\\
&+[\sum_{\substack{x\in A\\y\in B}}p(x,y)
g_{n}(A\Delta\{x,y\},t-u;{\bf s-u})g_{n}(B,t-u;{\bf s-u})]\notag\\
&+[\sum_{\substack{x\in B\\y\in A}}p(x,y)
g_{n}(A,t-u;{\bf s-u})g_{n}(B\Delta\{x,y\},t-u;{\bf s-u})] \Big\}.\notag
\end{align}
Define
\[
G_n(t)=\sup_{0\le s^1<\dots<s^n\le t}\sup_{\substack{
C\cap D=\emptyset\\ C, D\in \P^*}}
g_n(C\cup D,t;{\bf s})-g_n(C,t;{\bf s})g_n(D,t;{\bf s}),
\]
also set $F_n(t)=\sup\{G_n(s): 0\le s\le t\}$. 
Now, the key observation of Andjel, in~\cite{an} p. 720, is
that for $x\in A$ and $y\in B$ (so that both $x,y\not= 0$)
\[
g_n(A\cup B,t;{\bf s})\le F_n(t)+\frac{1}{2}\left(
g_n(A,t;{\bf s})g_{n}(B\Delta\{x,y\},t;{\bf s})+
g_n(A\Delta\{x,y\},t;{\bf s})g_{n}(B,t;{\bf s})\right).
\]
Thus, using the induction hypothesis ($F_{n-1}=0$) and
the symmetry of $p(.,.)$, we obtain
\begin{align}
g_n&(A\cup B,t;{\bf s})-g_n(A,t;{\bf s})g_n(B,t;{\bf s})\notag\\
&\le \int_0^{s^1}\!\! e^{-\l u}\left[ 
\sum_{x,y\in A}\!\! p(x,y)+\!\!\!\!\sum_{x,y\in B} \!\!p(x,y)+
\!\!\!\!\sum_{\substack{x\in A\\y\not\in A\cup B}}\!\!p(x,y)+
\!\!\!\!\sum_{\substack{x\in B\\y\not\in A\cup B}}\!\!p(x,y)+
2\sum_{\substack{x\in A\\y\in B}}p(x,y)\right] F_n(t-u)du\notag\\
&= \int_0^{s^1} e^{-\l u} \l F_n(t-u)du\le F_n(t) 
\int_0^{t} \l e^{-\l u}du\notag.
\end{align}
Thus, by taking the supremum over $A,B\in \P^*$ with $A\cap B=\emptyset$, we
obtain $G_n(t)\le F_n(t)\int_0^{t} \l \exp(-\l u)du$. Finally, this implies
that $F_n(t)=0$, and the proof is completed.
\qed
\subsection{Upper bound $g(\{i\},t)\le C_d \PP(H_i<\infty)$.}
\label{upperbound}
We first use the classical representation of
the trajectories $\{\L(\{i\},t),t>0\}$ as sequences
of jump times $\{\tau_i, i\in\NN\}$, which are independent
exponential variables of parameter $2d$, associated with the paths
of a simple symmetric random walk killed at the origin, say $\{\L_i,i\in\NN\}$.
The processes $\{\L_i,i\in\NN\}$ and $\{\tau_j, j\in\NN\}$ are independent.
We use the notation $E^y$ and $P^y$ to denote average over paths of
$\{\L_i,i\in\NN\}$ starting on $y$.
Let $T_0=\inf\{n>0:\ \L_n\in \N_0\}$ with $\N_0=\{y:y\sim 0\}$.
When $T_0<\infty$, let $T_1$ be the first return time to $\N_0$,
whereas when $T_0=\infty$, set $T_1=\infty$. Then, the sequence
of successive entrance times in $\N_0$, $\{T_2,T_3,\dots\}$ is defined
inductively. The number of return to $\N_0$ is called $R$, i.e. if
$T_i=\infty$ but $T_{i-1}<\infty$ then $R=i$. Our walk on $\ZZ^d_*$
being transient, we have $R<\infty$, a.s.  We note also, that by symmetry,
for any $y,y'\in \N_0$
\[
P^{y}(T_0<\infty)=P^{y'}(T_0<\infty)
\quad{\rm and}\quad\forall k\in \NN,\quad 
P^{y}(R=k)=P^{y'}(R=k).
\]
For convenience, we call $P^{\N_0}(T_0<\infty):=P^{y}(T_0<\infty)$
and $P^{\N_0}(R=k)=P^{y}(R=k)$ for $y\in \N_0$. Now,
\ba{eq1.10}
g(\{i\},t)&\le& E^{i}\Big[ 1_{\{T_0<\infty\}}\exp(2\sum_{i=0}^R\tau_{T_i})\Big]
=\sum_{k=0}^{\infty} E^{i}\Big[ 1_{\{T_0<\infty\}}1_{\{R=k\}}
\exp(2\sum_{i=0}^k\tau_{T_i})\Big]\cr
&=&P^{i}(T_0<\infty)\sum_{k=0}^{\infty} P^{\N_0}(R=k)(E[e^{2\tau_1}])^k.
\ea
where we used the Markov property and induction. Now, by the same
arguments
\[
P^{\N_0}(R=k)\le \left(P^{\N_0}(T_0<\infty)\right)^k.
\]
On the other hand, the evaluation of $E[\exp(2\tau_1)]$ is easy
\be{eq1.11}
E[e^{2\tau_1}]=\int_0^{\infty}e^{2t}2de^{-2dt}dt=\frac{d}{d-1}.
\ee
Thus, with (\ref{eq1.10}) and (\ref{eq1.11}), our upper bound follows
easily as soon as
\be{eq1.12}
P^{\N_0}(T_0<\infty)<\frac{d-1}{d}.
\ee
We want to formulate (\ref{eq1.12}) in terms of hitting probabilities for
the standard random walk, say $\{S_n,n\ge 0\}$.
We will denote the averages over the standard walk with a tilde.
Let $\kappa=\inf\{n>0: S_n=0\}$, and note first that
\[
P^{\N_0}(T_0<\infty)=\tilde P^{\N_0}(T_0<\infty, \kappa>T_0).
\]
By conditionning on the first move, we obtain
\[
\tilde P^{\N_0}(T_0<\infty)
=\frac{1}{2d}+\tilde P^{\N_0}(T_0<\infty, \kappa>T_0).
\]
Thus, (\ref{eq1.12}) is equivalent to $\tilde P^{\N_0}(T_0<\infty)<(2d-1)/2d$.
We recall that $R$ is the number of return to $\N_0$ for a walk
starting on $\N_0$. We can set $S_0=0$ but count what happens only after
two steps
\be{eq1.13}
R=\sum_{n=2}^{\infty} 1_{\{S_n\in \N_0\}}=\sum_{n\ge 1} 1_{\{T_n<\infty\}}.
\ee
So that 
\[
\tilde E[R]=\sum_{n=2}^{\infty}\tilde P(S_n\in \N_0)=
\frac{\tilde P^{\N_0}(T_0<\infty)}{1-\tilde P^{\N_0}(T_0<\infty)}.
\]
Thus, (\ref{eq1.12}) reads $\tilde E[R]<2d-1$. Now, we note that for $n>0$
\ba{eq1.14}
\tilde P(S_{n+1}=0)&=&\tilde P(S_{n+1}=0,S_n=y,\ y\in \N_0)\cr
&=&\sum_{y\in \N_0}
\tilde P(S_{n+1}=0|S_n=y)\tilde P(S_n=y)=\frac{1}{2d}\tilde P(S_n\in \N_0).
\ea
Thus,
\[
\sum_{n=2}^{\infty}\tilde P(S_n\in \N_0)=2d
\sum_{n=3}^{\infty}\tilde P(S_n=0).
\]
In dimension 3, it has been established (see \cite{durrett} page 170,
exercise 2.7) that
\[
1+\sum_{n=1}^{\infty}\tilde P(S_n=0)=
({\sqrt 6}/32\pi^3)\Gamma(1/24)\Gamma(5/24)\Gamma(6/24)\Gamma(11/24)=1.516...
\]
Thus, in dimension 3, $\tilde E[R]\le 6(0.52)<5$ and
our condition (\ref{eq1.12}) holds. 
Thus, in dimension 3, $\tilde E[R]\le 6(0.52)<5$ and our condition
(\ref{eq1.12}) holds. We conclude that (\ref{eq1.12}) holds for any dimension
larger or equal to 3 because the right hand side of $\tilde E[R]<2d-1$
increases and the average number of visits to 0 decreases with dimensions.
As we have not found a reference of this latter fact, we present a
short proof due to Andjel. The number of visits to $0$ is a geometric random 
variable of parameter $P_0(T_0< \infty)$, so it suffices to show 
monotonicity for this quantity. Let $i$ be any neighbor of 0; by symmetry 
$a(d):=P_0(T_0< \infty)=P_i(T_0<\infty)$, for dimension $d\ge3$.
Project the $d$-dimensional walk on a
hyperplane orthogonal to $i'$, a neighbor of the origin different from $i$. The
projected walk on $Z^{d-1}$ has transition probabilities $1/2d$ to go to each
of its $2(d-1)$ neighbors and $1/d$ not to move. It is clear that
$a(d)$ is not larger than the probability of visiting
the origin starting at (the projection of) $i$ for the projected walk.
This latter probability is actually equal to $a(d-1)$.
Indeed, the projected process goes along the same trajectories as
the $(d-1)$-dimensional standard walk and the waiting times at each point 
are geometric random variables with parameter $(d-1)/d$: thus, if the
trajectory of the standard walk is such that $\{T<\infty\}$, then the
same holds for the projected walk and vice-versa.

\section{Apriori Bounds.}
\label{apriori}
\subsection{Proof of Corollary~\ref{lem1}}
(i) The proof proceeds along the same lines 
as the proof of Theorem 3c) of~\cite{ad}, once the measure
$\nu_{\a(.)}$ is defined.
We set $\a_0=0$, and for $i\not= 0$, let $\a_i$ be defined by
\[
\frac{\a_i}{1-\a_i}\frac{1-\rho}{\rho}=\frac{1}{1+C_d\PP(H_i<\infty)},
\]
where the constant $C_d$ is that of (\ref{keypoint}). Note that $0<\a_i<\rho$.
Now, in the proof of Theorem 3c) of~\cite{ad}, we showed that
$\nu_{\a(.)}\ll\nur$ and
$d\nu_{\a(.)}/d\nur \in L^p(\nur)$ for $p>1$ as soon as
\[
\sum_{i\in \ZZ^d}(1-\frac{\a_i}{\rho})^2<\infty,
\quad\text{or equivalently in }d\ge 3\quad
\sum_{i\in \ZZ^d} \PP(H_i<\infty)^2<\infty,
\]
which holds as soon as $d\ge 5$.

We rewrite (\ref{keypoint}) on $\{\eta:\ \eta_i=0\}$ with $i\not=0$,
denoting by $\s_i$ the action of spin flip at site $i$ (i.e.
$(\sigma_i\eta)_k=\eta_k$ if $k\not= i$ and $(\sigma_i\eta)_i=1-\eta_i$)
\be{kp-bis}
\s_if_t\ge \frac{\a_i}{1-\a_i}\frac{1-\rho}{\rho} f_t\quad
\text{or equivalently}\quad \s_i (f_t) \frac{d\nu_{\a(.)}}{d\nur}\ge
\s_i(\frac{d\nu_{\a(.)}}{d\nur}) f_t.
\ee
Now, on $\A$, we form $\v=dT_t(\nur)/d\nu_{\a(.)}$ and we note that
$\v$ is increasing. Indeed, if $i\not= 0$ and $\eta_i=0$,
then (\ref{kp-bis}) is nothing but $\s_i \v\ge \v$.
Now, as a product
measure $\nu_{\a(.)}$ satisfies FKG. Thus, for $\psi$ increasing
\[
\int \psi dT_t(\nur)=
\int \psi(\frac{dT_t(\nur)}{d\nu_{\a(.)}})d\nu_{\a(.)}\ge
\int \psi d\nu_{\a(.)}\int(\frac{dT_t(\nur)}{d\nu_{\a(.)}})d\nu_{\a(.)}=
\int \psi d\nu_{\a(.)}.
\]
Thus, $\nu_{\a(.)}\prec T_t(\nur)$. The fact that $T_t(\nur)\prec \nur$ comes
from the fact that $f_t$ is decreasing and $\nur$ satisfies FKG.
Now, (ii) of Corollary~\ref{lem1} follows as in the
proof of Theorem 3c) of~\cite{ad}. Using that $f_t$ and $d\nu_{\a(.)}/d\nur$
are decreasing, for $i\ge 1$ and $j\ge 0$,
\[
\int\!\! f_t^i \left(\frac{d\nu_{\a(.)}}{d\nu_{\rho}}\right)^j
d\nu_{\rho} = \int\!\! f_t^{i-1}
\left(\frac{d\nu_{\a(.)}}{d\nu_{\rho}}\right)^j\!\! dT_t(\nur)\leq
\int\!\! f_t^{i-1}
\left(\frac{d\nu_{\a(.)}}{d\nu_{\rho}}\right)^j \!\! d\nu_{\a(.)} =
\int\!\! f_t^{i-1} \left(\frac{d\nu_{\a(.)}}{d\nu_{\rho}}\right)^{j+1}
\!\! d\nu_{\rho} .
\]
Thus, we obtain by induction, for each $n \geq 1$
\be{ineq.ad}
\int f_t^n d\nu_{\rho} \leq \int
\left(\frac{d\nu_{\a(.)}}{d\nu_{\rho}}\right)^n d\nu_{\rho}.
\ee
Since the right hand side of (\ref{ineq.ad}) is bounded for $d \geq 5$,
the corollary follows.

\subsection{Proof of Lemma~\ref{lem.gene}}
In \cite{ad} section 4, we have that $t\mapsto R(t):=e^{-\lambda t}/
P_{\nu}(\tau>t)$ is increasing and $R(0)>0$. Suppose $(i)$ and
let $\lim_{t\to\infty} R(t)=R<\infty$. Now, $(ii)$ follows from 
\be{eq2.gene}
\int f_t^2d\nu=
\int \frac{\bar S_t(1_{\A^c})\bar S_t(1_{\A^c})}{P_{\nu}(\tau>t)^2}
d\nu=\int \frac{\bar S_{2t}(1_{\A^c})}{P_{\nu}(\tau>t)^2}
d\nu={P_{\nu}(\tau>2t)\over P_{\nu}(\tau>t)^2}=\frac{R(2t)}
{R(t)^2}\le \frac{R}{R(0)^2}.
\ee

Conversely, Let $\mu$ be a limit point of $\{1/t\int_0^t T_s(\nu)ds\}$
along the subsequence $\{t_n\}$ in weak-$L^2(d\nu)$. By Theorem 1
of \cite{ad}, $\mu$ is a quasi-stationary measure with
$P_{\mu}(\tau>t)=\exp(-\lambda t)$. Thus, $(i)$ follows from
\be{eq3.gene}
\lim_{t\to\infty}R(t)=\lim_{n \ra +\infty} {1\over t_n} \int_0^{t_n}
{P_{\mu}(\tau > s)\over P_{\nu_{\rho}}(\tau > s)} ds
=\lim_{n \ra +\infty} {1\over t_n} \int_0^{t_n} \int f_s
{d\mu\over d\nu_{\rho}} d\nu_{\rho} ds=
\int \left( {d\mu\over d\nu_{\rho}} \right)^2 d\nu_{\rho}<\infty.
\ee
Now, if either $(i)$ or $(ii)$ holds,
by remark 3 of \cite{ad}, we have that the Yaglom limit exists.
Now, to show that $\int(f_t-f)^2d\nu$ converges to 0, we only
need to show that $\int f_t^2d\nu$ converges to $\int f^2d\nu$.
It is easy to see from \eqref{eq2.gene} and \eqref{eq3.gene}, that
$\lim \int f_t^2d\nu=\lim R(t)=\int (d\mu/d\nu)^2 d\nu$, and the
proof of Lemma~\ref{lem.gene} is concluded.

\section{Proof of Theorem~\ref{the2} (i): Uniqueness}
\label{unique}
\subsection{Positivity of $f$.}
We first show that $\nur$-a.s., $f>0$ on $\{\eta_0=0\}$. We introduce a
symmetric simple exclusion process on $\ZZ^d\bs\{0\}$: 
there is no site 0, and its adjacent bonds 
are suppressed. Let $\nur^*$ be the product Bernoulli measure
on $\ZZ^d\bs\{0\}$ of density $\rho$. 
Let $\{S^*_t,\ t\ge 0\}$ be the Markov semi-group of this process.
It is known that the process $\{S^*_t,\ t\ge 0\}$,
with initial measure $\nur^*$, is reversible and ergodic:
indeed, by Theorem 1.44 on page 377 of~\cite{l1}, $\nur^*$ is
an extremal invariant measure and by Theorem B.52, on page 23 of~\cite{l2},
all extremal invariant measures are ergodic. In other words, if for any $t>0$,
$S^*_tg =g$, $\nur^*$-a.s., then $g$ is constant $\nur^*$-a.s.

Let $\B:=\{\eta: \eta_0=0, f(\eta)=0\}$, and note that $\nur$-a.s.
\[
\bar S_tf. 1_\B=\exp(-\lambda t) f.1_\B=0. 
\]
Now, because $f\ge 0$, we have that for any $\e>0$, 
$\e 1_{\{f>\e\}} \le f$. Thus, for $\eta\in \B$, $\nur$-a.s.
\[
\forall \e>0,\quad P^{\eta}(\{f(\eta_t)>\e\}, \tau>t)=0
\Longrightarrow P^{\eta}(\{f(\eta_t)>0\}, \tau>t)=0.
\] 
Now, if the bonds linking 0 to its neighbors are not marked, 
up to time $t$, then
$\bar S_t$ acts like $S^*_t$. So, if $\tau_0$ is the first time a Poisson mark
appears in one of these bonds, for $\eta\in \B$, $\nur$-a.s.
\[
P^{\eta}( f(\eta_t)>0, \tau_0>t)=0,
\]
because $\tau_0\le \tau$. Now, $\tau_0$ is independent from the Poisson
processes on bonds of $\ZZ^d\bs\{0\}$, so
\[
P^{\eta}( f(\eta_t)>0, \tau_0>t)=S^*_t 1_{\{f(.)>0\}} P(\tau_0>t)=0.
\]
Now, for any $t>0$, $P(\tau_0>t)>0$. Thus, for $\eta\in \B$, 
$\nur^*$-a.s.
\[
S^*_t1_{\{f(.)>0\}}(\eta)=0.
\] 
In other words, for any $t>0$, we have $\nur^*$-a.s.
\[
S^*_t1_\B \ge 1_\B.
\]
Now, $\nur^*$ is reversible for $S^*_t$, so both expressions 
have the same mean, and we conclude that $\nur^*$-a.s., 
for any $t>0$, $S^*_t1_\B=1_\B$.
By the ergodicity of $\nur^*$, we conclude that $1_\B$ is $\nur^*$
constant, so that necessarily $\nur(\B)=0$.
\subsection{One eigenvalue with a positive eigenvector.}
Suppose that $f,f'\in \H_\A$ are the densities of
two quasi-stationary measures. There are
two real numbers $\lambda(\rho)$ and $\lambda'$ such that
$f,f'$ satisfy in an $L^2(\nur)$ sense
\be{eq2.1}
\bar S_tf=e^{-\lambda(\rho) t} f,
\quad{\rm and}\quad
\bar S_tf'=e^{-\lambda' t} f'.
\ee
First, we show that $\lambda(\rho)=\lambda'$. We assume that
$d\mur=fd\nur$ corresponds to the Yaglom limit. Thus,
\be{eq2.2}
\lim_{t\to\infty} \int f'f_td\nur=\int f'fd\nur
,\quad{\rm and}\quad
\lim_{t\to\infty} \int ff_td\nur=\int f^2 d\nur.
\ee
However, as $f>0$ on $\{\eta_0=0\}$, $\nur$-a.s.,
\be{eq2.3}
{e^{-\lambda't}\int f'd\nur\over P_{\nur}(\tau>t)}
={\int 1_{A^c}\bar S_tf'd\nur\over P_{\nur}(\tau>t)}=
\int f'f_td\nur \longrightarrow\int f'fd\nur>0.
\ee
Similarly,
\be{eq2.4}
\int ff_td\nur={e^{-\lambda(\rho) t}\int fd\nur\over P_{\nur}(\tau>t)}
\longrightarrow\int f^2d\nur.
\ee
Thus, $\lambda(\rho)=\lambda'$. 
\subsection{Dual expansion.}
We expand $f$ on the countable basis of $L^2(\nur)$, say $\{H_A,\ A\in\P\}$,
where $\P$ is the collection of finite subsets of $\ZZ^d$ and 
\be{eq2.7}
H_{\emptyset}=1,\quad{\rm and}\quad
H_A(\eta)=\prod_{i\in A} {(\eta_i-\rho)\over \sqrt{\rho
(1-\rho)}}.
\ee
Thus, there are real numbers $\{C_A,\ A\in\P\}$ with
\be{eq2.8}
f=\sum_{A\in\P} C_A H_A,\quad{\rm and}\quad
\int f^2d\nur=\sum_{A\in\P} C_A^2.
\ee
The constraint that $f\in \H_\A$, i.e. $\eta_0f(\eta)=0$, is equivalent to
\be{eq2.9}
\forall A\not\ni 0,\quad \int H_A \eta_0 fd\nur=0.
\ee
Thus,
\be{eq2.10}
\sqrt{\rho(1-\rho)}\int H_{A\cup\{0\}}fd\nur+\rho\int H_A fd\nur=0.
\ee
So, for $A\not\ni 0$,
\be{eq2.11}
\sqrt{\rho(1-\rho)}C_{A\cup\{0\}}+\rho C_A=0.
\ee
We define for all $A\in \P^*$
\be{eq2.12}
\psi(A)=\left(-\sqrt{{1-\rho\over \rho}}\right)^{|A|}C_A.
\ee
Now, condition (\ref{eq2.11}) reads $\psi(A\cup\{0\})=\psi(A)$ for
$A\not\ni 0$. Now, we express $(-Lf,f)$ in terms of the $\{C_A\}$.
\be{eq2.13}
(-Lf,f)=\sum_{A\in \P}\sum_{B\sim A}(C_B-C_A)^2,
\ee
where $B\sim A$ if there is $i\in A\bs B$ and $j\in B\bs A$ with 
$A\bigtriangleup B=\{i,j\}$. We replace the $C_A$'s by the $\psi(A)$'s
and distinguish 0 to eliminate (\ref{eq2.11}) 
\ba{eq2.14}
(-Lf,f)&=& \sum_{A\not\ni 0}
\left( \sum_{B\not\ni 0,B\sim A}\g^{|A|}(\psi(B)-\psi(A))^2+
\sum_{B\ni 0,B\sim A}\g^{|A|}(\psi(B\bs\{0\})-\psi(A))^2\right)\cr
&+&\sum_{A\ni 0}\left(
\sum_{B\ni 0,B\sim A}\g^{|A|}(\psi(B\bs\{0\})-\psi(A\bs\{0\}))^2+
\sum_{B\not\ni 0,B\sim A}\g^{|A|}(\psi(B)-\psi(A\bs\{0\}))^2\right)\cr
&=&\sum_{A\not\ni 0}\Big[
\sum_{B\in \N_A}\g^{|A|}(1+\g)(\psi(B)-\psi(A))^2+
\sum_{B\in \N_A^-}\g^{|A|}(\psi(B)-\psi(A))^2\cr
&& +\sum_{B\in \N_A^+}\g^{|A|+1}(\psi(B)-\psi(A))^2\Big],
\ea
where $\g=(1-\rho)/\rho$, $B\in \N_A$ means $B\sim A$ and $B\not\ni 0$,
$B\in \N_A^+$ means $B\sim A\cup\{0\}$ and $B\not\ni 0$, and
$B\in \N_A^-$ means $B\cup\{0\}\sim A$ and $B\not\ni 0$. 
With this rewriting, $(-\bar Lf,f)$ can be thought of as
the Dirichlet form, $\E(\psi,\psi)$,
of a dynamics with finitely many particles, with
creation and annihilation at site 0,
with respect to a measure $m(A)=\g^{|A|}$.
The advantages of this rewriting are threefold: (i)
the constraint (\ref{eq2.11}) has vanished, (ii) the new dynamics
is clearly irreducible
and (iii) in studying the minimizers of $\E(\psi,\psi)$,
we can assume the $\{\psi(A)\}$ to be nonnegatives. Indeed, note that
$\E(\psi,\psi)\ge \E(|\psi|,|\psi|)$ and equality holds only
if $\psi(A)\ge 0$ for each $A\in \P^*$.
Also, we rewrite the $L^2(\nur)$ norm of $f$ in terms of the
$\{\psi(A), A\in \P^*\}$
\be{eq2.15}
|f|^2=\sum_{A\not\ni 0} C_A^2+C_{A\cup\{0\}}^2=(1+\g)\sum_{A\not\ni 0}
\psi(A)^2 \g^{|A|}=(1+\gamma)||\psi||_m^2.
\ee
Now, $\bar Lf+\lambda f=0$ implies that for $A\in \P^*$
\be{eq2.16}
(1+\g)\!\!\sum_{B\in \N_A}(\psi(B)-\psi(A))+
\g\!\!\sum_{B\in \N_A^+}(\psi(B)-\psi(A))+
\!\!\sum_{B\in \N_A^-}(\psi(B)-\psi(A))=-\lambda \psi(A).
\ee
Thus, $\psi(A)>0$ for all $A\in \P^*$. 
Now, let $\phi=\psi^2$ and note
that for any $A$ and $B\in \P^*$, the functional
$\phi\mapsto (\sqrt{\phi(B)}-\sqrt{\phi(A)})^2$ is convex. Assume
that $\psi$ and $\psi'$ are two positive normalized minimizers, and let
$\phi$ and $\phi'$ be their respective squares. Then, for any
$\lambda\in [0,1]$
\[
\psi_{\lambda}:=\sqrt{\lambda \phi+(1-\lambda)\phi'}\ {\rm has}\ 
||\psi_{\lambda}||_{m}=1 \quad{\rm and}\quad
\E(\psi_{\lambda},\psi_{\lambda})\le \lambda\E(\psi,\psi)+
(1-\lambda)\E(\psi',\psi').
\]
Thus, the convex inequality is an equality, so that for
any $A\in \P^*$ and any $B\in \N_A\cup \N_A^-\cup \N_A^+$
\[
(\psi_{\lambda}(B)-\psi_{\lambda}(A))^2=
\lambda(\psi(B)-\psi(A))^2+ (1-\lambda)(\psi'(B)-\psi'(A))^2,
\]
which implies, after expanding, that $\phi(A)\phi'(B)=\phi(B)\phi'(A)$. Now,
as $\phi(A)>0$ and the dynamics is irreducible, we conclude that
$\phi\equiv\phi'$ and so are the positive square-roots $\psi\equiv \psi'$.

\section{Proof of Theorem~\ref{the2} (ii) and (iii)}
\label{density}
\subsection{Proof of Theorem~\ref{the2} (ii): Density at infinity}
The facts that for any $t>0$, $\nu_{\a(.)}\prec T_t(\nur)\prec \nur$ with
$\a(i)\to \rho$ when $||i||\to\infty$ implies that for any $A\in\P$
\[
\int \prod_{j\in A} \eta_{j+i}d\mur(\eta)
\quad\substack{||i||\to\infty \\ \longrightarrow} \quad
\rho^{|A|},
\]
for $\prod_{j\in A} \eta_{j+i}$ is an increasing function. Now any local
function $\v$ can be written as a linear combination of local
increasing functions, and the property follows by linearity.

\subsection{Proof of Theorem~\ref{the2} (iii): Basin of Attraction}
We show that for any measure $\nu$, any subsequence of the Ces\`aro limit
of $\{T_t(\nu)\}$ contains a further subsequence converging
to a quasi-stationary measure, say $\mu$.
When the density of $\nu$ with respect to $\nur$, say $\phi$, is
continuous, we show that $\mu=\mur$.
As in the proof of existence of a quasi-stationary
measure (see \cite{ad} Lemma1), we establish first the existence, for
any $s>0$, of the following limit
\be{eq4.1}
\lim_{t\to\infty} {P_{\nu}(\tau>t+s)\over P_{\nu}(\tau>t)}=\exp(-
\lambda(\rho) s).
\ee
Indeed, recall that 
\be{eq4.2}
P_{\nu}(\tau>t)=\int \bar S_t(1_{A^c})1_{A^c}\phi d\nur=
\int \bar S_t(1_{A^c}\phi)d\nur,
\ee
so that by the existence of the Yaglom limit
\be{eq4.7}
\lim_{t\to\infty}{P_{\nu}(\tau>t)\over P_{\nur}(\tau>t)}=
\int \phi d\mur>0.
\ee
Thus, 
\be{eq4.4}
\lim_{t\to\infty} {P_{\nu}(\tau>t+s)\over P_{\nu}(\tau>t)}=
\lim_{t\to\infty} {P_{\nur}(\tau>t+s)\over P_{\nur}(\tau>t)}=
\exp(-\lambda(\rho) s).
\ee
By the weak$^*$ compactness, for any sequence
$\{t_n\}$, there is a further subsequence
(still named $\{t_n\}$ for convenience), and
a probability measure $\mu$, such that for any continuous 
function $\v$
\be{eq4.5}
{1\over t_n}\int_0^{t_n} (T_t(\nu),\v)dt\longrightarrow \int \v d\mu.
\ee
The same argument
as in the proof of Theorem 1 of \cite{ad} implies that $\mu$ is
quasi-stationary. Now, we establish apriori estimates
\ba{eq4.6}
\int\left({dT_t(\nu)\over d\nur}\right)^2d\nur&=&
\int\left({E^{\eta}[\phi(\eta_t)1_{\{\tau>t\}}]\over P_{\nu}(\tau>t)}
\right)^2d\nur\cr
&\le&|\phi|_{\infty}^2\int\left({P^{\eta}(\tau>t)
\over P_{\nu}(\tau>t)}\right)^2d\nur
=|\phi|_{\infty}^2\int f_t^2d\nur \left({P_{\nur}(\tau>t)
\over P_{\nu}(\tau>t)}\right)^2.
\ea
This quantity is bounded by Corollary~\ref{lem1} and (\ref{eq4.7}).
By standard arguments, this implies that $\mu\in L^2(\nur)$
which by the uniqueness result establish that $\mu=\mur$, so that
the Ces\`aro limit exists and is $\mur$.
\vspace{0,5cm}

\noindent{\bf Acknowledgements}. A.A. thanks for their warmth, 
the faculty and staff of
the mathematics department of S\~ao Paulo University, where part
of this research was conducted. Also, we would like
to thank Fabienne Castell and Enrique Andjel for suggestions
and Paolo Dai Pra for pointing out an error in the first version.

\end{document}